\documentclass{article}
\usepackage{amsmath,amssymb}

\newcommand{\C}{{\mathbb C}}

\author{Mark van Hoeij}

\title{\sc An algorithm for computing
        the Weierstrass normal form of hyperelliptic curves}

\begin{document}
\date{\empty}
\hyphenation{}
\maketitle
\begin{abstract}
An algorithm is given to compute a normal form for hyperelliptic curves.
The elliptic case has been treated in~\cite{issac}. In this paper the
hyperelliptic case is treated.
\end{abstract}

\section{Introduction}
\label{s1}
The problem we consider is the same as in~\cite{issac},
but then for the hyperelliptic case instead of the elliptic case.
Let $C$ be an algebraic curve given by an irreducible polynomial
$f \in \C[x,y]$. The function field of $C$ is the field of fractions
of $\C[x,y]/(f)$, and is denoted as $\C(x,y)$.
We want to decide if $C$ is hyperelliptic or not.
If $C$ is hyperelliptic, we want to find a Weierstrass normal form
for $C$ (i.e. an equation $y_0^2 = F(x_0)$ for some squarefree
polynomial $F(x_0) \in \C[x_0]$).
We want to find a birational map from $C$ to this normal form,
and its inverse. This means expressing $x,y$ in terms of $x_0,y_0$
and vice versa. Let $g$ be the genus of $C$.

The problems to be solved are the following
\begin{itemize}
\item Decide if $C$ is hyperelliptic or not. By definition, $C$
is hyperelliptic if $g>1$ and there exists $x_0 \in \C(x,y)$
such that $\C(x_0)$ is a subfield of $\C(x,y)$ with index 2.
If $C$ is hyperelliptic then the field $K := \C(x_0)$ is unique.
\item Find a generator $x_0$ of $K$.
\item Find an element $y_0$ of $\C(x,y)$, not in $K$.
Take for example $y_0=x$  (interchange $x,y$ if it turns out
that $x \in K$).
Compute the minimal polynomial of $y_0=x$ over $K$.
This is done by writing $x_0$ as a quotient $v_1/v_2$ with
$v_1,v_2 \in \C[x,y]$ and computing the resultant
${\rm Res}_y(v_2 \tilde{x}_0 - v_1, f)$ where $\tilde{x}_0$
is the name of a variable that plays the role of $x_0$.
Then add an element of $K$
to $y_0$ to change its minimal polynomial into the form
$y_0^2-F(x_0)$ for some rational function $F(x_0) \in K$. Then multiply 
$y_0$ by an element of $K$ to turn $F_0$ into a squarefree polynomial.
\item At this point the functions $x_0,y_0$ are expressed as
elements of $\C(x,y)$ and their algebraic relation $y_0^2-F(x_0)$
is in Weierstrass form. Now view $x_0,y_0$ as variables $\tilde{x}_0,
\tilde{y}_0$, and express $x,y$ as rational functions in
these variables. This can be done in the same way as in~\cite{issac}.
\end{itemize}

\section{Decide if $C$ is hyperelliptic}
If $f \in \C[x,y]$ is irreducible, then we can decide if $C$
is elliptic by computing the genus $g$. The curve is elliptic
if and only if $g=1$.
For the hyperelliptic case more work needs to be done.
An algorithm for this is implemented in Maple 6.
To view the it, issue the following commands in Maple 6 or higher \\
\verb+interface(verboseproc=2): op(`algcurves/is_hyperelliptic`);+ \\
and the following code appears
\begin{verbatim}
proc(f, x, y)
local a, b, i, j, g, z;
option remember, `Copyright (c) 1999 Waterloo Maple Inc. All rights\
 reserved. Author: M. van Hoeij`;
    i := degree(f, {x, y});
    j := min(degree(f, x), degree(f, y), i - ldegree(f, {x, y}));
    if i < 4 or j < 2 then false
    elif j = 2 then evalb(1 < `algcurves/genus`(f, x, y))
    else
        b := `algcurves/differentials`(f, x, y, skip_dx);
        g := nops(b);
        j := 1/2*(i - 1)*(i - 2) - g;
        if g = 2 then true
        elif g < 2 or j < 2 or j < 4 and
        {seq(z[3], z = `algcurves/singularities`(f, x, y))} = {1}
        then false
        else
            z := {seq(a[i], i = 1 .. g)};
            j := {coeffs(expand(numer(subs(RootOf(f, y) = y, evala(
                Expand(
                subs(y = RootOf(f, y), add(a[i]*b[i], i = 1 .. g)^2))))))
                , z)};
            z := {seq(a[i], i = 1 .. nops(j))};
            z := nops(j) - nops(indets(subs(`solve/linear`({coeffs(
                expand(add(a[i]*j[i], i = 1 .. nops(j))), {x, y})}, z)
                , z)) intersect z);
            if z < 2*g - 1 then error "should not happen"
            else evalb(z = 2*g - 1)
            end if
        end if
    end if
end proc
\end{verbatim}
We will not explain all steps, but the main steps are
the following.
\begin{itemize}
\item Let $V$ be the vector space of holomorphic differentials.
An algorithm described in~\cite{differentials} is used
to compute a basis $b_1,\ldots,b_g$ for $V$.
\item If $g<2$ then the curve is not hyperelliptic and the algorithm
stops.
\item Let $W$ be the span of all $b_i b_j$, $1 \leq i \leq j \leq g$.
\item If $C$ is hyperelliptic, and if $t$ is any generator of $K$, i.e.
$K = \C(t)$, then it is known that
\begin{equation}
	t^i {\rm d}t, \ i=0,\ldots,g-1
\label{bas}
\end{equation}
is a basis of $V$. Then
$t^i {\rm d}^2t$, $i=0,\ldots,2g-2$ must be a basis of $W$.
So if $C$ is hyperelliptic then the dimension of $W$ must be $2g-1$.
It is not hard to show that if $C$ is not hyperelliptic and $g>1$
then ${\rm dim}(W) > 2g-1$. \\
In the algorithm, $z = {\rm dim}(W)$. Testing if $z=2g-1$
decides if $C$ is hyperelliptic or not.
\end{itemize}

\section{Compute a function of degree 2}
Once we have a generator $x_0$ for $K$, it is not difficult
to find $y_0$, see section~\ref{s1}. This is implemented in Maple
in the procedure \verb+`algcurves/genus2`+  (to view the code type
similar commands as in the previous section).
After that, the remaining steps are the same as in~\cite{issac}.

So the main remaining problem to be solved is to compute some generator
$x_0$ of $K$. We have a basis $b_1,\ldots,b_g$ of $V$.
We want to find two particular elements of $V$, whose quotient
is a generator of $K$. For example, if $t$ is any generator of $K$,
then ${\rm d}t$, $t {\rm d}t$
will do, and so will $t^{g-2} {\rm d}t$, $t^{g-1} {\rm d}t$.
We will describe how to find a subspace
$V_2 \subseteq V$ of the form $\C t^{g-2} {\rm d}t + \C t^{g-1} {\rm d}t$
for some generator $t$ of $K$.
If $v_1,v_2$ is any basis of $V_2$, then we can
take $x_0 = v_1/v_2$ as our generator of $K$.

The idea for finding $V_2$ is the following. Take an arbitrary nonsingular 
point $P$ on the curve. Since $t$ may be any generator of $K$, we
may assume without loss of generality that $t(P) = 0$.
Then $t^{g-1} {\rm d}t$ has the highest valuation at $P$ of
any nonzero element of $V$, and $t^{g-2} {\rm d}t$ has the second
highest valuation. So we can find $V_2$ from $V$ by solving linear 
equations: As long as ${\rm dim}(V) > 2$, let $m$ be the smallest 
valuation at $P$ of any element of $V$. Then compute the subspace of $V$ 
of those elements that have valuation $>m$, and replace $V$ by that 
subspace. This way the dimension decreases, and after a number of steps 
$V$ will equal $V_2$. Then take $x_0=v_1/v_2$ where $v_1,v_2$ is any
basis of $V_2$.
This is implemented in Maple in the procedure
\verb+`algcurves/index2subfield`+.

\section{An example}
\[ f=y^9+3x^2y^6+3x^4y^3+x^6+y^2 \]
This curve has genus 3 and is hyperelliptic.
The result of the Maple commands: \\
\verb+with(algcurves): Weierstrassform(f,x,y,x0,y0);+\\ is
\begin{eqnarray*}
[{{y_0}}^{2}+2-7\,{x_0}+35\,{{x_0}}^{4}-21\,{{x_0}}^{5}+7
\,{{x_0}}^{6}-{{x_0}}^{7}+21\,{{x_0}}^{2}-35\,{{x_0}}^{3}, \\[5pt]
{\frac {y \left( {y}^{3}+{x}^{2}+1 \right) }{{y}^{6}+2\,{x}^{2}{y}^{3}
+{x}^{4}+y}},
-{\frac {x \left( {y}^{3}+{x}^{2} \right) }{y}},- \left( 
-1+{x_0} \right) {y_0},1-3\,{x_0}+3\,{{x_0}}^{2}-{{x_0}
}^{3}]
\end{eqnarray*}
The first expression is the Weierstrassform $y_0^2 - F(x_0)$.
The second expression is a generator for the field $K$, i.e. a function
of degree 2 in $\C(x,y)$. Under the birational map this corresponds
to $x_0$. The third expression is the image of $y_0$ under this
birational map. The 4'th and 5'th expression are $x,y$ expressed
in terms of $x_0,y_0$.

\end{document}